\documentclass[11pt]{amsart} 
\usepackage{amssymb,amsmath,latexsym} 

\newtheorem{theorem}{Theorem} 
\newtheorem{lemma}[theorem]{Lemma} 
\newtheorem{corollary}[theorem]{Corollary} 
\newtheorem*{algorithm}{Algorithm}

\title[Partial-fractions method for integral linear systems]{The partial-fractions method for counting solutions to integral linear systems} 

\author{Matthias Beck}
\address{Department of Mathematics, San Francisco State University, 1600 Holloway Avenue, San Francisco, CA 94132, U.S.A.}
\email{beck@math.sfsu.edu}

\begin{document}
\setlength{\parindent}{0pt} 

\newcommand\const{\operatorname{const}} 
\newcommand\Hil{\operatorname{Hil}} 
\newcommand\rank{\operatorname{rank}} 
\renewcommand\th{^{\text{th}}} 
\def\A{{\bf A}} 
\def\b{{\bf b}}
\def\c{{\bf c}}
\def\n{{\bf n}}
\def\r{{\bf r}}
\def\v{{\bf v}}
\def\x{{\bf x}}
\def\z{{\bf z}}
\def\C{\mathbb{C}}
\def\R{\mathbb{R}}
\def\Z{\mathbb{Z}}
\def\P{\mathcal{P}}
\def\Q{\mathcal{Q}}
\def\T{\mathcal{T}}
\def\L{L}

\abstract 
We present a new tool to compute the number $\phi_\A (\b)$ of integer solutions to the linear system
$$
\x \geq 0 \ , \qquad \A \, \x = \b \ , 
$$
where the coefficients of $\A$ and $\b$ are integral. $\phi_\A (\b)$ is often described as a \emph{vector partition function}. Our methods use partial fraction expansions of Euler's generating function for $\phi_\A (\b)$. A special class of vector partition functions are Ehrhart (quasi-)polynomials counting integer points in dilated polytopes.
\endabstract 

\dedicatory{Dedicated to Lou Billera on the occasion of his sixtieth birthday.}

\subjclass[2000]{\emph{Primary} 05A15, 52C07; \emph{Secondary} 52C45. 
} 

\keywords{Vector partition function, partial fractions expansion, quasi-polynomial, lattice-point counting, rational convex polytope, Ehrhart theory. 
}
\thanks{This paper was written while the author was at the Mathematical Sciences Research Institute in Berkeley and at the Max-Planck-Institut f\"ur Mathematik in Bonn. The author thanks both institutes for their hospitality. Thanks also to Ira Gessel and Sinai Robins for fruitful discussions, and to the referees for many helpful comments.}
\thanks{Appeared in \emph{Discrete \& Computational Geometry} {\bf 32} (2004), 437--446 (special issue in honor of Louis Billera).}
\date{February 29, 2004}

\maketitle  

\setlength{\parskip}{0.4cm}
\bibliographystyle{amsplain}


\section{Euler's generating function}

We are interested in computing the number of integer solutions of the linear system 
$$
\x \in \R_{\geq 0}^d \ , \qquad \A \, \x = \b \ , 
$$
where $\A$ is a nonnegative $(m \times d)$-integral matrix and $\b \in \Z^m$. 
We think of $\A$ as fixed and study the number of solutions $\phi_\A(\b)$ as a function of $\b$. 
(Strictly speaking, this function should only be defined for those $\b$ which lie in the nonnegative linear span of the columns of $\A$.) 

The function $\phi_\A(\b)$, often called a \emph{vector partition function}, appears in a wealth of mathematical areas and beyond: Number Theory (partitions), Discrete Geometry (polyhedra), Commutative Algebra (Hilbert series), Algebraic Geometry (toric varieties), Representation Theory (tensor product multiplicities), Optimization (integer programming), as well as applications to Chemistry, Biology, Physics, Computer Science, and Economics. 

Denote the columns of $\A$ by $\c_1, \dots, \c_d$. The following lemma goes back to at least Euler \cite{euler}: 
\begin{lemma}[Euler]
\label{eulerlemma} 
$\phi_\A(\b)$ equals the coefficient of $\z^\b := z_1^{b_1} \cdots z_m^{b_m}$ of the function 
$$ 
f(\z) = \frac 1 { \left( 1 - \z^{ \c_1 } \right) \cdots \left( 1 - \z^{ \c_d } \right) } 
$$ 
expanded as a power series centered at $\z = 0$. 
\end{lemma} 

\begin{proof}
Expand each factor of the right-hand side into a geometric series.
\end{proof}

Equivalently, the coefficient of $\z^\b$ in $f(\z)$ equals the constant term in $\frac{ f(\z) }{ \z^\b }$, denoted by $\const \frac{ f(\z) }{ \z^\b }$. So Euler's Lemma can be conveniently stated as
$$ 
\phi_\A(\b) = \const \frac 1 { \left( 1 - \z^{ \c_1 } \right) \cdots \left( 1 - \z^{ \c_d } \right) \z^\b } \ .
$$
In a series of articles \cite{tetrahed,bdr,beckpixton,polygons}, we used complex integration of $\frac{ f(\z) }{ \z^\b }$ to compute $\phi_\A(\b)$ for special cases of $\A$. Similar techniques were applied in \cite{lasserrezeron}. Here we expand $\frac{ f(\z) }{ \z^\b }$ into partial fractions to compute its constant term, and hence $\phi_\A(\b)$. This work constitutes in a sense a refinement of the complex-integration methods, with the advantage that it is more flexible and---more importantly---applicable for \emph{any} integral linear system. 



\section{Vector partition functions}

The nature of the counting function $\phi_\A (\b)$ is given by the following theorem.
A \emph{quasi-polynomial} is a finite sum $Q(\b)=\sum_\n c_\n (\b) \, \b^\n$ with coefficients $c_\n$ that are functions of $\b$ which are periodic in every component of $\b$. The \emph{degree} of $Q$ is  the degree of the largest power $\b^\n$ appearing in $Q$. 
A matrix is \emph{unimodular} if every square submatrix has determinant $\pm 1$. 

\begin{theorem}[Sturmfels \cite{sturmfelsvectorpartition}]
\label{sturmthm}
The function $\phi_\A (\b)$ is a piecewise-defined quasi-polynomial in $\b$ of degree $d-\rank(\A)$. 
The regions of $\R^m$ in which $\phi_\A (\b)$ is a single quasi-polynomial are polyhedral, that is, they are defined by linear constraints. 
If $\A$ is unimodular then $\phi_\A$ is a piecewise-defined polynomial. 
\end{theorem}

The unimodular case of this theorem is due to Dahmen and Micchelli \cite[Corollary 3.1]{dahmenmicchelli}.

The computation of both $\phi_\A (\b)$ and the chamber complex consisting of the regions of quasi-polynomiality give rise to challenging problems. The most powerful technique for computing $\phi_\A (\b)$ which we are aware of is due to Brion and Vergne \cite{brionvergne}. (The methods described here are much more elementary.) 
The chamber complex is still much of a mystery. A promising approach can be found in \cite{szenesvergne}.

We finish this section with a \emph{reciprocity theorem}. Let $\phi_\A^\circ (\b)$ count the integer solutions of 
$$
\x > 0 \ , \qquad \A \, \x = \b \ .
$$
Both $\phi_\A$ and $\phi_\A^\circ$ are quasi-polynomials, and can hence be algebraically defined for arguments which are not integer vectors in the positive span of $\A$. The following identity shows the close relationship between the two functions.

\begin{theorem}[\cite{evencloser}]
\label{multirec}
The quasi-polynomials $\phi_\A$ and $\phi_\A^\circ$ satisfy
$$
\phi_\A(-\b) = (-1)^{ d - \rank \A } \phi_\A^\circ (\b) \ .
$$
\end{theorem}

This identity gives rise to a symmetry property of $\phi_\A(\b)$. Let $r_k$ denote the sum of the entries in the $k\th$ row of $\A$, and let $\r = (r_1,\dots,r_m)$. Then the integer solutions of 
$$
\x > 0 \ , \qquad \A \, \x = \b
$$
are in bijection (via $x_k \mapsto x_k-1$) with the integer solutions of 
$$
\x \geq 0 \ , \qquad \A \, \x = \b - \r \ ,
$$
and hence $\phi_\A^\circ (\b) = \phi_\A (\b-\r)$. This yields:
\begin{corollary}
\label{reccor}
The quasi-polynomial $\phi_\A$ satisfies
$$
\phi_\A(\b) = (-1)^{ d - \rank \A } \phi_\A (-\b-\r) \ . 
$$
\end{corollary}


\section{The partial-fractions method}

This section describes the idea behind our computations. Recall that our goal is to derive
$$ 
\phi_\A(\b) = \const \frac 1 { \left( 1 - \z^{ \c_1 } \right) \cdots \left( 1 - \z^{ \c_d } \right) \z^\b } \ .
$$
We start by expanding
$$ 
\frac 1 { \left( 1 - \z^{ \c_1 } \right) \cdots \left( 1 - \z^{ \c_d } \right) \z^\b } \ .
$$
into partial fractions in one of the components of $\z$, say $z_1$: 
$$ 
\frac 1 { \left( 1 - \z^{ \c_1 } \right) \cdots \left( 1 - \z^{ \c_d } \right) \z^\b } 
= \frac 1 { z_2^{b_2} \cdots z_m^{b_m} } \sum_{k=1}^d \frac{ A_k (\z,b_1) }{ 1-\z^{\c_k} } + \sum_{j=1}^{b_1} \frac{ B_j (\z) }{ z_1^j } \ .
$$
Here $A_k$ and $B_j$ are polynomials in $z_1$, rational functions in $z_2, \dots, z_m$, and exponential in $b_1$. (We are tacitly assuming that there are no multiple poles (besides $\z=0$), which generally holds unless $m=1$, a case which can be handled easily.) The two sums on the right-hand side correspond to the analytic and the meromorphic part with respect to $z_1=0$. The latter does not contribute to the $z_1$-constant term, whence
\begin{align}
\phi_\A (\b)
&= \const_{z_2, \dots, z_m} \left( \frac 1 { z_2^{b_2} \cdots z_m^{b_m} } \, \const_{z_1} \left( \sum_{k=1}^d \frac{ A_k (\z,b_1) }{ 1-\z^{\c_k} } \right) \right) \nonumber \\
&= \const \left( \frac 1 { z_2^{b_2} \cdots z_m^{b_m} } \sum_{k=1}^d A_k (0,z_2,\dots,z_m,b_1) \right) \ . \label{mainidea}
\end{align}
The effect of one partial fraction expansion is to eliminate one of the variables of the 
generating function, at the cost of replacing one rational function by a sum of such. 
The following idea is now evident. 

\begin{algorithm}
Apply (\ref{mainidea}) repeatedly to eliminate $z_1$, then $z_2$, etc., up to $z_{m-1}$. 
\end{algorithm}

The constant term of the remaining rational functions in the one variable $z_m$ can be computed with the methods introduced in \cite{bdr}. This implies, in particular, that any vector partition function is a quasi-polynomial whose nontrivial ingredients are \emph{Fourier-Dedekind sums} \cite{bdr}. The constant-term computation for the last variable also explains the quasi-polynomial character of $\phi_\A$, since in this last step roots of unity will appear, with the components of $\b$ as exponents.

We note that this algorithm computes $\phi_\A(\b)$ \emph{as a function of $\b$}, that is, it allows symbolic computation. 
Secondly, it is not very hard to deduce Sturmfels's Theorem \ref{sturmthm} from this algorithm. 
What might be more important, however, is the fact that the constraints which define the regions of quasi-polynomiality of $\phi_\A$ are obtained ``on the go" as one computes $\phi_\A$: When expanding into partial fractions, one has to check where the poles of a rational function are. The components of $\b$ will appear (linearly) in the exponents of these rational functions, and hence one will automatically have to split the computation into cases which give rise to different regions in the chamber complex of quasi-polynomiality.

While we hope that our algorithm gives a computational tool, in particular, for obtaining the chamber complex defining the regions of quasi-polynomiality of a vector partition function, it is not clear to us why the chambers obtained from the partial fractions analysis coincide with Sturmfels's predicted chambers.

Our algorithm is best illustrated by going through an actual example.


\section{An illustrating example}
\label{examplesec}

Let 
$\A = \left( \begin{array}{cccc} 1 & 2 & 1 & 0 \\ 1 & 1 & 0 & 1 \end{array} \right)$, and write 
$\b = (a,b)$, so $\phi_\A(a,b)$ counts the integer solutions of 
$$
x_1, x_2, x_3, x_4 \geq 0 \ , \qquad 
\begin{array}{rcl}
x_1 + 2x_2 + x_3 & = & a \\
x_1 + x_2 + x_4  & = & b \ .
\end{array}
$$

By Euler's Lemma \ref{eulerlemma}, 
$$
\phi_\A (a,b) = \const \frac 1 { (1-zw) (1-z^2 w) (1-z) (1-w) z^a w^b } \ . 
$$
We first expand into partial fractions with respect to $w$:
$$
\frac 1 { (1-zw) (1-z^2 w) (1-w) w^b } = - \frac{ \frac{ z^{b+1} }{ (1-z)^2 } }{ 1-zw } + \frac{ \frac{ z^{ 2b+3 } }{ (1-z) (1-z^2) } }{ 1 - z^2 w } + \frac{ \frac{ 1 }{ (1-z) (1-z^2) } }{ 1-w } + \sum_{k=1}^b \frac{ \dots }{ w^k } \ .
$$
Taking constant terms gives
\begin{align}
\phi_\A (a,b) &= \const_z \left( \frac 1 { (1-z) z^a } \, \const_w \left( \frac 1 { (1-zw) (1-z^2 w) (1-w) w^b } \right) \right) \nonumber \\
&= \const \left( \frac 1 { (1-z) z^a } \left( - \frac{ z^{b+1} }{ (1-z)^2 } + \frac{ z^{ 2b+3 } }{ (1-z) (1-z^2) } + \frac{ 1 }{ (1-z) (1-z^2) } \right) \right) \nonumber \\
&= \const \left( - \frac{ z^{b-a+1} }{ (1-z)^3 } + \frac{ z^{ 2b-a+3 } }{ (1-z)^2 (1-z^2) } + \frac{ 1 }{ (1-z)^2 (1-z^2) z^a } \right) \label{exparfrac}
\end{align}

(At this point, we could interpret each of the three constant terms as counting integer solutions to new linear systems. For example, the last term gives the number of integer points $(x,y) \geq 0$ satisfying $2x+y \leq a$. To keep a general flavor, we continue with our general algorithm.)

We compute the constant term of each of the three terms separately. For starters,
$$
\const \frac{ z^{b-a+1} }{ (1-z)^3 } = 0
$$
if $b-a+1 > 0$, equivalently (since $a$ and $b$ are integers) $b \geq a$. If $b<a$, we use
$$
\frac 1 {(1-z)^3} = \sum_{k \geq 0} \binom {k+2} 2 z^k \ ,
$$
which gives 
$$
\const \frac{ 1 }{ (1-z)^3 z^{a-b-1} } = \binom {a-b+1} 2 = \frac {(a-b)^2} 2 + \frac {a-b} 2 \ . 
$$
Incidentally, this identity is true not only for $b<a$, but also for $b=a$ and $b=a+1$, because the right-hand side vanishes then. This suggests that the regions of quasi-polynomiality \emph{overlap}, as was proved by Szenes and Vergne \cite{szenesvergne}. Hence for the first constant term, we obtain
$$
\const \frac{ 1 }{ (1-z)^3 z^{a-b-1} } =
\begin{cases}
0 & \text{ if } b \geq a, \\
\frac {(a-b)^2} 2 + \frac {a-b} 2 & \text{ if } b \leq a+1.
\end{cases}
$$

For the second term in (\ref{exparfrac}), 
$$
\const \frac{ z^{ 2b-a+3 } }{ (1-z)^2 (1-z^2) } = 0
$$
if $2b-a+2 \geq 0$. If $a \geq 2b+3$ we expand into partial fractions again: 
\begin{align*} 
\const \frac{ z^{ 2b-a+3 } }{ (1-z)^2 (1-z^2) }
&= \const \left( \frac{ 1/2 }{ (1-z)^3 } + \frac{ \frac {a-2b-3} 2 + \frac 1 4 }{ (1-z)^2 } + \frac{ \frac {(a-2b-3)^2} 4 + \frac {a-2b-3} 2 + \frac 1 8 }{ 1-z } + \frac{ (-1)^{a+1}/8 }{ 1+z } \right) \\
&= \frac {(a-2b)^2} 4 + \frac {2b-a} 2 + \frac { 1 + (-1)^{a+1} } 8 \ .
\end{align*}
Similar to the first constant-term computation, this identity is also valid for $a=2b+2, 2b+1, 2b$, whence
$$
\const \frac{ z^{ 2b-a+3 } }{ (1-z)^2 (1-z^2) } =
\begin{cases}
0 & \text{ if } a \leq 2b+2, \\
\frac {(a-2b)^2} 4 + \frac {2b-a} 2 + \frac { 1 + (-1)^{a+1} } 8 & \text{ if } a \geq 2b.
\end{cases}
$$

The computation for the last term in (\ref{exparfrac}) is almost identical. (Note that this term always contributes.)
\begin{align*} 
\const \frac{ 1 }{ (1-z)^2 (1-z^2) z^a}
&= \const \left( \frac{ 1/2 }{ (1-z)^3 } + \frac{ \frac a 2 + \frac 1 4 }{ (1-z)^2 } + \frac{ \frac {a^2} 4 + \frac a 2 + \frac 1 8 }{ 1-z } + \frac{ (-1)^a/8 }{ 1+z } \right) \\
&= \frac {a^2} 4 + a + \frac { 7 + (-1)^a } 8 \ .
\end{align*}

Summing up all terms in (\ref{exparfrac}) gives finally
\begin{align*} 
\phi_\A (a,b) &= \begin{cases}
\frac {a^2} 4 + a + \frac { 7 + (-1)^a } 8 & \text{ if } a \leq b , \\
- \frac {(a-b)^2} 2 - \frac {a-b} 2 + \frac {a^2} 4 + a + \frac { 7 + (-1)^a } 8 & \text{ if } \frac a 2 - 1 \leq b \leq a+1 , \\
- \frac {(a-b)^2} 2 - \frac {a-b} 2 + \frac {(a-2b)^2} 4 + \frac {2b-a} 2 + \frac { 1 + (-1)^{a+1} } 8 + \frac {a^2} 4 + a + \frac { 7 + (-1)^a } 8 & \text{ if } b \leq \frac a 2 \\
\end{cases} \\
&= \begin{cases}
\frac {a^2} 4 + a + \frac { 7 + (-1)^a } 8 & \text{ if } a \leq b , \\
ab - \frac {a^2} 4 - \frac {b^2} 2 + \frac {a+b} 2 + \frac { 7 + (-1)^a } 8 & \text{ if } \frac a 2 - 1 \leq b \leq a+1 , \\
\frac {b^2} 2 + \frac {3b} 2 + 1 & \text{ if } b \leq \frac a 2 . \\
\end{cases}
\end{align*}

It is a fun exercise to show that $\phi_\A (a,b) = \phi_\A (-a-4,-b-3)$, as promised by Corollary \ref{reccor}.

We pause for a moment to introduce the geometry behind our linear system. By thinking of $x_3$ and $x_4$ as slack variables, we can see that
$$
x_1, x_2, x_3, x_4 \geq 0 \ , \qquad 
\begin{array}{rcl}
x_1 + 2x_2 + x_3 & = & a \\
x_1 + x_2 + x_4  & = & b \ .
\end{array}
$$
is equivalent to 
$$
x_1,x_2 \geq 0 \ , \qquad 
\begin{array}{rcl}
x_1 + 2x_2 & \leq & a \\
x_1 +  x_2 & \leq & b \ .
\end{array}
$$
Depending on the relationship between $a$ and $b$, the geometric figure described by this linear system is a quadrilateral or triangle. 

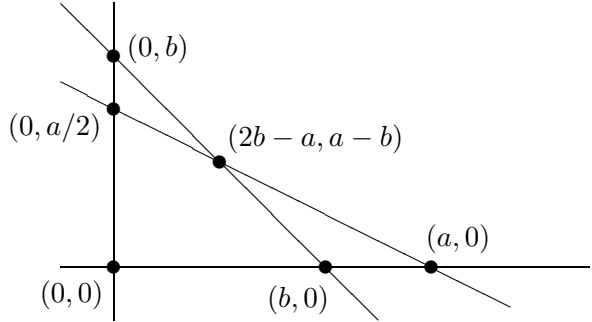
\begin{center} 
\begin{picture}(200,120)
\put(5,20){\line(1,0){200}}
\put(25,0){\line(0,1){120}}
\put(5,120){\line(1,-1){120}}
\put(5,90){\line(2,-1){170}}

\put(25,100){\circle*{5}}
\put(30,100){$(0,b)$}
\put(25,80){\circle*{5}}
\put(-15,70){$(0,a/2)$}
\put(105,20){\circle*{5}}
\put(83,5){$(b,0)$}
\put(145,20){\circle*{5}}
\put(143,28){$(a,0)$}
\put(65,60){\circle*{5}}
\put(68,66){$(2b-a,a-b)$}
\put(25,20){\circle*{5}}
\put(-3,7){$(0,0)$}
\end{picture} 
\end{center} 

The inequalities which define the different cases for $\phi_\A (a,b)$ determine the chamber complex in the parameter space, and in this example also the combinatorial type of the polytope: Except for boundary cases, in the first ($a \leq b$) and last ($b \leq \frac a 2$) case, it is a triangle, and in the second case ($a > b > \frac a 2$) a quadrilateral. In fact, the geometry of the triangle in the first case only depends on $a$, which is reflected in $\phi_\A (a,b)$, and similarly the last case shows only dependency on $b$. In this last case, the triangle actually has integer vertices (the linear system is unimodular), whence we get a polynomial. In the first two cases the vertices of the polytope are half-integral, which is reflected in the period-2 quasi-polynomials. The following picture illustrates the three overlapping chambers.

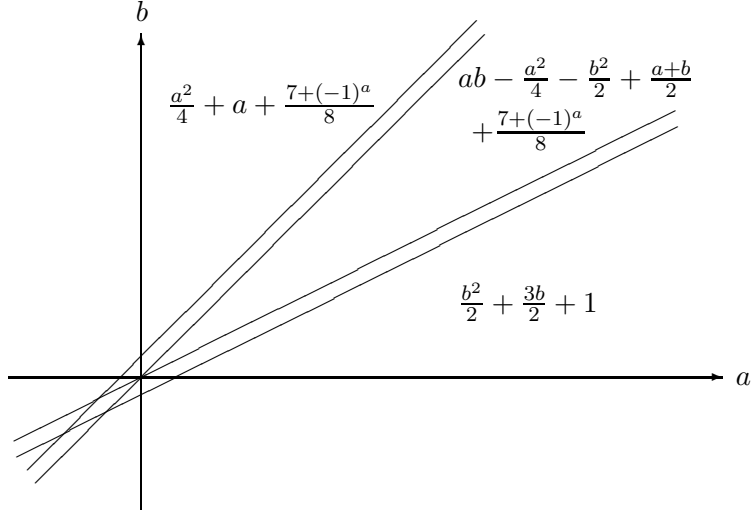
\begin{center} 
\begin{picture}(280,190)
\put(0,50){\vector(1,0){270}}
\put(275,47){$a$}

\put(50,0){\vector(0,1){180}}
\put(48, 185){$b$}

\put(10,10){\line(1,1){170}} 
\put(7,15){\line(1,1){170}} 
\put(3,20){\line(2,1){250}} 
\put(2,26){\line(2,1){250}} 

\put(60,150){$\frac {a^2} 4 + a + \frac { 7 + (-1)^a } 8$}
\put(170,160){$ab - \frac {a^2} 4 - \frac {b^2} 2 + \frac {a+b} 2$} 
\put(175,140){$+ \frac { 7 + (-1)^a } 8$}
\put(170,75){$\frac {b^2} 2 + \frac {3b} 2 + 1$}
\end{picture} 
\end{center} 


\section{Ehrhart quasi-polynomials}

A \emph{convex polytope} $\P$ in $\R^d$ is the convex hull of finitely many points in $\R^d$. Alternatively (and this correspondence is nontrivial \cite{ziegler}), one can define $\P$ as the bounded intersection of affine halfspaces. 
A polytope is \emph{rational} if all of its vertices have rational coordinates. (A \emph{vertex} of $\P$ is a point $\v \in \P$ for which there is a hyperplane $H$ such that $\{ v \} = \P \cap H$.) 
We denote by $\P^\circ$ the relative interior of $\P$. 
For a positive integer $t$, let $ L_\P (t) $ denote the number of integer points (``lattice points") in the dilated polytope $ t \P = \{ tx : x \in \P  \} $. The fundamental result about the structure of $L_\P$ is as follows. 
\begin{theorem}[Ehrhart \cite{ehrhart2}]
\label{quasi}
If $\P$ is a convex rational polytope, then the functions $ L_\P (t) $ and $ L_{\P^\circ} (t) $ are quasi-polynomials in $t$ whose degree is the dimension of $\P$. 
If $\P$ has integer vertices, then $ L_\P $ and $ L_{\P^\circ} $ are polynomials. 
\end{theorem} 
Ehrhart conjectured and partially proved the following \emph{reciprocity law}, which was proved by Macdonald \cite{macdonald}. 
\begin{theorem}[Ehrhart-Macdonald]
\label{rec}
The quasi-polynomials $L_\P$ and $L_{\P^\circ}$ satisfy
$$
  L_\P (-t) = (-1)^{\dim \P} L_{\P^\circ} (t) \ .
$$
\end{theorem} 

The computation of Ehrhart quasi-polynomials is only slightly easier than that of $\phi_\A$. Recent work includes \cite{barvinokehrhart,brionvergne,cappellshaneson,chendcg02,latte,diazrobins,guillemin,kantorkhovanskii,khovanskipukhlikov,pommersheim}. 

Suppose the convex polytope $\P \subset \R^d$ is given by an intersection of halfspaces, that is, 
$$
\P = \left\{ \x \in \R^d : \ \A \, \x \leq \b \right\} \ ,
$$
for some $(m \times d)$-matrix $\A$ and $m$-dimensional vector $\b$. We may convert these inequalities into equalities by introducing slack variables. If $\P$ has rational vertices, we can choose $\A$ and $\b$ in such a way that all their entries are integers, without loss of generality nonnegative ones. In summary, we may assume that a convex rational polytope $\P$ is given by 
$$
\P = \left\{ \x \in \R_{ \geq 0 }^d : \ \A \, \x = \b \right\} \ ,
$$
where $\A \in M_{m \times d} (\Z_{ \ge 0 } )$ and $\b \in \Z_{ \ge 0 }^m$. (If we are interested in counting the integer points in $\P$, we may assume that $\P$ is in the 
nonnegative orthant, i.e., the points in $\P$ have nonnegative coordinates, as translation by an integer vector does not change the lattice-point count.) 

The connection to vector partition functions is now evident. Since 
$
t \P = \left\{ \x \in \R_{ \geq 0 }^d : \ \A \, \x = t \, \b \right\} 
$,
we obtain $\L_\P (t) = \phi_\A (t \b)$ as a special evaluation of $\phi_\A$. Note that $t \b, t=1,2, \dots$ lie in the same chamber of quasi-polynomiality of $\phi_\A$. Ehrhart's Theorem \ref{quasi} is therefore a special case of Sturmfels's Theorem \ref{sturmthm}, and Theorem \ref{rec} is a special case of Theorem \ref{multirec}. 

As an example, the quadrilateral $\Q$ described by
$$
x,y \geq 0 \ , \qquad 
\begin{array}{rcl}
x + 2y & \leq & 5 \\
x + y  & \leq & 4 
\end{array}
$$
(a special case of the polygons appearing in Section \ref{examplesec}) with vertices $(0,0),(4,0),(3,1),(0,5/2)$ has the Ehrhart quasi-polynomial
$$
\L_\Q (t) = \phi_\A (5t,4t) = \frac {23} 4 \, t^2 + \frac 9 2 \, t + \frac { 7 + (-1)^t } 8 \ .
$$

Ehrhart quasi-polynomials are easier to compute, since one does not need to derive/know the chamber complex of quasi-polynomiality of $\phi_\A$. Our algorithm, naturally, works just as well for $\L_\P$. 


\section{Concluding remarks}

Many open problems and questions remain untouched. 
Most importantly, we hope that our ideas will give rise to practical implementations of computing vector partition functions  (including the chamber complex defining the regions of quasi-polynomiality) in their various disguises in form of a computer software. Again, the partial-fraction method is more general and more flexible than algorithms using complex residues, for example in \cite{beckpixton}.

Not unrelated is the question of computational complexity. It is known that the rational generating function of the Ehrhart quasi-polynomial of a $d$-polytope can be computed in polynomial time if $d$ is fixed \cite{barvinokalgorithm}. We have not analyzed the complexity of our algorithm. This seems interesting in the light that our algorithm depends on the number of contraints, geometrically corresponding to the facets (codimension-1 faces) of the polytope. Barvinok's algorithm, in contrast, depends on the number of vertices of the polytope. 

Another venue which should be explored is the following: After each step in our algorithm, one could reinterpret the constant term of each summand as counting integer solutions of a linear system. It might be interesting to simplify each of these linear systems (as far as this is possible without changing the number of integer solutions), and then proceed with the algorithm.


\bibliographystyle{amsplain}

\begin{thebibliography}{10}

\bibitem{barvinokehrhart}
Alexander~I. Barvinok, \emph{Computing the {E}hrhart polynomial of a convex
  lattice polytope}, Discrete Comput. Geom. \textbf{12} (1994), no.~1, 35--48.

\bibitem{barvinokalgorithm}
\bysame, \emph{A polynomial time algorithm for counting integral points in
  polyhedra when the dimension is fixed}, Math.~Oper.~Res. \textbf{19} (1994),
  769--779.

\bibitem{tetrahed}
Matthias Beck, \emph{Counting lattice points by means of the residue theorem},
  Ramanujan J. \textbf{4} (2000), no.~3, 299--310.

\bibitem{evencloser}
\bysame, \emph{Multidimensional {E}hrhart reciprocity}, J. Combin. Theory Ser.
  A \textbf{97} (2002), no.~1, 187--194.

\bibitem{bdr}
Matthias Beck, Ricardo Diaz, and Sinai Robins, \emph{The {F}robenius problem,
  rational polytopes, and {F}ourier-{D}edekind sums}, J. Number Theory
  \textbf{96} (2002), no.~1, 1--21.

\bibitem{beckpixton}
Matthias Beck and Dennis Pixton, \emph{The {E}hrhart polynomial of the
  {B}irkhoff polytope}, Discrete Comput. Geom. \textbf{30} (2003), no.~4,
  623--637.

\bibitem{polygons}
Matthias Beck and Sinai Robins, \emph{Explicit and efficient formulas for the
  lattice point count in rational polygons using {D}edekind-{R}ademacher sums},
  Discrete Comput. Geom. \textbf{27} (2002), no.~4, 443--459.

\bibitem{brionvergne}
Michel Brion and Mich{\`e}le Vergne, \emph{Residue formulae, vector partition
  functions and lattice points in rational polytopes}, J. Amer. Math. Soc.
  \textbf{10} (1997), no.~4, 797--833.

\bibitem{cappellshaneson}
Sylvain~E. Cappell and Julius~L. Shaneson, \emph{Euler-{M}aclaurin expansions
  for lattices above dimension one}, C. R. Acad. Sci. Paris S\'er. I Math.
  \textbf{321} (1995), no.~7, 885--890.

\bibitem{chendcg02}
Beifang Chen, \emph{Lattice points, {D}edekind sums, and {E}hrhart polynomials
  of lattice polyhedra}, Discrete Comput. Geom. \textbf{28} (2002), no.~2,
  175--199.

\bibitem{dahmenmicchelli}
Wolfgang Dahmen and Charles~A. Micchelli, \emph{The number of solutions to
  linear {D}iophantine equations and multivariate splines}, Trans. Amer. Math.
  Soc. \textbf{308} (1988), no.~2, 509--532.

\bibitem{latte}
J.~DeLoera, R.~Hemmecke, J.~Tauzer, and R.~Yoshida, \emph{Effective lattice
  point counting in rational convex polytopes}, Preprint, 2003.

\bibitem{diazrobins}
Ricardo Diaz and Sinai Robins, \emph{The {E}hrhart polynomial of a lattice
  polytope}, Ann. of Math. (2) \textbf{145} (1997), no.~3, 503--518.

\bibitem{ehrhart2}
Eug{\`e}ne Ehrhart, \emph{Sur un probl\`eme de g\'eom\'etrie diophantienne
  lin\'eaire. {I}{I}. {S}yst\`emes diophantiens lin\'eaires}, J. Reine Angew.
  Math. \textbf{227} (1967), 25--49.

\bibitem{euler}
Leonhard Euler, \emph{De partitione numerorum in partes tam numero quam specie
  dates}, Leonhardi Euleri Opera Omnia, Ser.~1 Vol.~III (F.~Rudio, ed.),
  Teubner, Leipzig, 1917, pp.~132--147.

\bibitem{guillemin}
Victor Guillemin, \emph{Riemann-{R}och for toric orbifolds}, J. Differential
  Geom. \textbf{45} (1997), no.~1, 53--73.

\bibitem{kantorkhovanskii}
Jean-Michel Kantor and Askold Khovanskii, \emph{Une application du th\'eor\`eme
  de {R}iemann-{R}och combinatoire au polyn\^ome d'{E}hrhart des polytopes
  entiers de $\mathbf {R}\sp d$}, C. R. Acad. Sci. Paris S\'er. I Math.
  \textbf{317} (1993), no.~5, 501--507.

\bibitem{khovanskipukhlikov}
A.~G. Khovanski{\u\i} and A.~V. Pukhlikov, \emph{The {R}iemann-{R}och theorem
  for integrals and sums of quasipolynomials on virtual polytopes}, Algebra i
  Analiz \textbf{4} (1992), no.~4, 188--216.

\bibitem{lasserrezeron}
Jean~B. Lasserre and Eduardo~S. Zeron, \emph{On counting integral points in a
  convex rational polytope}, Math. Oper. Res. \textbf{28} (2003), no.~4,
  853--870.

\bibitem{macdonald}
I.~G. Macdonald, \emph{Polynomials associated with finite cell-complexes}, J.
  London Math. Soc. (2) \textbf{4} (1971), 181--192.

\bibitem{pommersheim}
James~E. Pommersheim, \emph{Toric varieties, lattice points and {D}edekind
  sums}, Math. Ann. \textbf{295} (1993), no.~1, 1--24.

\bibitem{sturmfelsvectorpartition}
Bernd Sturmfels, \emph{On vector partition functions}, J. Combin. Theory Ser. A
  \textbf{72} (1995), no.~2, 302--309.

\bibitem{szenesvergne}
Andr{\'a}s Szenes and Mich{\`e}le Vergne, \emph{Residue formulae for vector
  partitions and {E}uler-{M}ac{L}aurin sums}, Adv. in Appl. Math. \textbf{30}
  (2003), no.~1-2, 295--342, Formal power series and algebraic combinatorics
  (Scottsdale, AZ, 2001).

\bibitem{ziegler}
G{\"u}nter~M. Ziegler, \emph{Lectures on polytopes}, Springer-Verlag, New York,
  1995.

\end{thebibliography}

\def\cprime{$'$}
\providecommand{\bysame}{\leavevmode\hbox to3em{\hrulefill}\thinspace}
\providecommand{\MR}{\relax\ifhmode\unskip\space\fi MR }
\providecommand{\MRhref}[2]{%
  \href{http://www.ams.org/mathscinet-getitem?mr=#1}{#2}
}
\providecommand{\href}[2]{#2}

\setlength{\parskip}{0cm}

\end{document}